\documentclass[11pt]{amsart}
\usepackage{eurosym}
\usepackage{verbatim}
\usepackage{eucal,url,amssymb,stmaryrd,enumerate,amscd,verbatim}
\usepackage{amsmath}
\usepackage[pagebackref,colorlinks=true,linkcolor=blue,citecolor=blue]{hyperref}
\usepackage{amsfonts}
\usepackage{amsthm}
\usepackage[margin=1in]{geometry}

\allowdisplaybreaks
\numberwithin{equation}{section}
\newtheorem{thrm}{Theorem}[section]

\newtheorem{prop}[thrm]{Proposition}

\newcommand{\ddt}{\frac{d}{dt}}

 1

\def\gr{\nabla f}

\newcommand{\vol}{\, Vol_{\eta}}

\newcommand{\Cr}{\nabla}

\begin{document}

\begin{abstract}
 A proof of the monotonicity of an entropy like energy for the heat equation on a quaternionic contact and CR manifolds is proven.
\end{abstract}

\title[The monotonicity for the heat equation on a qc and CR manifolds]{The
monotonicity of an entropy like energy for the heat equation on a quaternionic contact and CR manifolds}
\date{\today }
\author{D. Vassilev}
\address[Dimiter Vassilev]{ Department of Mathematics and Statistics\\
University of New Mexico\\
Albuquerque, New Mexico, 87131-0001}
\email{vassilev@math.unm.edu}
\maketitle
\tableofcontents


\setcounter{tocdepth}{2}

\section{Introduction}

The purpose of this note is to show the monotonicity of the entropy type energy for the
heat equation on a compact quaternionic contact manifold inspired by the corresponding Riemannian fact related to Perelman's entropy formula for the heat equation on a static Riemannian manifold, see \cite{Ni04}. More recently a similar quantity was considered in the CR case \cite{CW10}. Our goal is to give a relatively simple  proof of the monotonicity, more in line with the Riemannian case, by resolving directly the difficulties arising in the sub-Riemannian setting. In Section \ref{s:CR} we include a proof of the result of \cite{CW10} in the CR case from our point of view.

To state the problem, let $M$ be a
quaternionic contact manifold, henceforth abbreviated to qc, and $u$ be a
smooth \emph{positive} solution to the quaternionic contact heat equation
\begin{equation}  \label{e:heat}
\frac{\partial}{\partial t}u=\Delta u.
\end{equation}
Hereafter, $\triangle u=tr^g(\nabla^2 u)$ is the negative sub-Laplacian with
the trace taken with respect to an orthonormal basis of the horizontal $4n$%
-dimensional space. Associated to such a solution are the (Nash like) entropy
\begin{equation}\label{e:N definition}
\mathcal{N}(t)=\int_M u\ln u \vol
\end{equation}
and entropy energy
functional
\begin{equation}\label{enf}
\mathcal{E}(t)=\int_M |\nabla f|^2 u \, Vol_{\eta},
\end{equation}
where, as usual, $f=-\ln u$ and $\, Vol_{\eta}$ is the naturally associated
volume form on $M$, see \eqref{e:volume form} and also \cite[Chapter 8]{IMV}. Exactly as in the Riemannian case, we have that the entropy is decreasing (i.e., non-increasing) because of the formula
\[
\ddt \mathcal{N}=-\mathcal{E}(t).
\]

Our goal is the computation of the second derivative of the entropy. In order to state the result we consider the Ricci type tensor
\begin{equation}  \label{e:Lichnerowicz tensor}
\mathcal{L}(X,X)\overset{def}{=}2 Sg(X,X)+\alpha_n T^0(X,X) +\beta_n
U(X,X)\geq 4g(X,X),
\end{equation}
where $X$ is any vector from the horizontal distribution, $\alpha_n=\frac {%
2(2n+3)}{2n+1}, \quad \beta_n=\frac {4(2n-1)(n+2)}{(2n+1)(n-1)}$, and $T^0$
and $U$ are certain invariant components of the torsion, see Subsection \ref{ss:prelim}.
In addition, following \cite{IPV2},
we define the $P-$form of a fixed smooth function $f$ on $M$ by the following equation
\begin{multline}  \label{e:P form}
P_f(X) =\sum_{b=1}^{4n}\nabla
^{3}f(X,e_{b},e_{b})+\sum_{t=1}^{3}\sum_{b=1}^{4n}\nabla
^{3}f(I_{t}X,e_{b},I_{t}e_{b}) \\
-4nSdf(X)+4nT^{0}(X,\nabla f)-\frac{8n(n-2)}{n-1}U(X,\nabla f),
\end{multline}
which in the case $n=1$ is defined by formally dropping the last term. The $P-$%
function of $f$ is the function $P_f(\nabla f)$. The $C-$operator of $M$ is the 4-th order differential operator
\begin{equation*}
f\mapsto Cf =-\nabla^* P_f=\sum_{a=1}^{4n}(\nabla_{e_a} P_f)\,(e_a).
\end{equation*}%
In many respects the $C-$operator plays a role similar to the Paneitz
operator in CR geometry. We say that the $P-$function of $f$ is non-negative
if 
\begin{equation*}
\int_M f\cdot Cf \, Vol_{\eta}= -\int_M P_f(\nabla f)\, Vol_{\eta} \geq 0.
\end{equation*}%
If the above holds for any $f\in \mathcal{C}^\infty_o\,(M)$ we say that the $%
C-$operator is {non-negative}, $\mathcal{C}\geq 0$.

We are ready to state our first result.
\begin{prop}
\label{p:energy ineq} Let $M$ be a compact QC manifold of dimension $4n+3$.
If $u=e^{-f}$ is a positive solution to heat equation \eqref{e:heat}, then
we have
\begin{equation*} 
\frac {2n+1}{4n}\mathcal{E}'(t)=-\int_M \left[ |(\nabla^2f)_{0}|^2+\frac {2n+1%
}{2}\mathcal{L}(\nabla f,\nabla f) +\frac {1}{16n}| \nabla f|^4 \right]u\,
Vol_{\eta} +\frac {3}{n}\int_M P_F(\nabla F)\, Vol_{\eta},
\end{equation*}
where $u=F^2$ ($f=-2\ln F$) and $(\nabla^2f)_{0}$ is the traceless part of horizontal Hessian of $f$.
\end{prop}

Several important properties of the C-operator were found in \cite{IPV2},
most notable of which is the fact that the $C-$operator is non-negative for $n>1$. In dimension seven, $n=1$, the condition of non-negativity of
the $C-$operator is non-trivial. However, \cite{IPV2} showed that
on a 7-dimensional compact qc-Einstein manifold with positive qc-scalar curvature
the $P-$function of an {eigenfunction} of the sub-Laplacian is
non-negative. In particular, this property holds on any 3-Sasakian manifold. Clearly, these facts together with Proposition \ref{p:energy ineq} imply the following theorem.

\begin{thrm}
\label{t:monotone energy} Let $M$ be a compact QC manifold of dimension $4n+3
$ of non-negative Ricci type tensor $\mathcal{L}(X,X) \geq 0$. In the case $%
n=1$ assume, in addition, that the $C-$operator is non-negative. If $u=e^{-f}
$ is a positive solution to heat equation \eqref{e:heat} then the energy is
monotone decreasing (i.e., non-increasing).
\end{thrm}

The proof of Proposition \ref{p:energy ineq} follows one of L. Ni's
arguments \cite{Ni04} in the Riemannian case, thus it relies on Bochner's
formula. More precisely, after Ni's initial step, in order to handle the extra terms in Bochner's formula, we will follow the
presentation of \cite{IV14}
where this was done for the qc Lichnerowicz type lower eigenvalue bound
under positive Ricci type tensor, see \cite{IPV1,IPV2} for the original result. In the qc case, similar
to the CR case, the Bochner formula has additional hard to control terms,
which include the $P$-function of $f$. In our case, since the integrals are
with respect to the measure $u\, Vol_{\eta}$, rather than $\vol$ as in the Lichnerowicz type estimate, some new
estimates are needed. The key is the following proposition which can
be considered as an estimates from above of the integral of the $P$-function
of $f$ with respect to the measure $u\, Vol_{\eta}$ when the $C-$operator is
non-negative.

\begin{prop}
\label{p:Paneitz estimate} Let $(M,\eta)$ be a compact QC manifold of
dimension $4n+3$. If $u=e^{-f}$ is a positive solution to heat equation %
\eqref{e:heat}, then with $f=-2\ln F$ we have the identity
\begin{equation}  \label{e:Paneitz estimate}
\int_M P_f(\nabla f) u\, Vol_{\eta}=\frac 14\int_M |\nabla f|^4 u\,
Vol_{\eta} +4\int_M P_F(\nabla F)\, Vol_{\eta}.
\end{equation}
\end{prop}

In the last section of the paper we apply the same method in the case of a strictly pseudoconvex pseudohermitian manifold and prove the following Proposition.

\begin{prop}\label{p:CR energy ineq} Let $M$ be a compact strictly pseudoconvex pseudohermitian CR manifold of dimension $2n+1$.
If $u=e^{-f}$ is a positive solution to the heat equation \eqref{e:heat}, then
we have
\begin{equation*} 
\frac {n+1}{2n}\mathcal{E}'(t)=-\int_M \left[ |(\nabla^2f)_{0}|^2+\frac {2n+1%
}{2}\mathcal{L}(\nabla f,\nabla f) +\frac {1}{8n}| \nabla f|^4 \right]u\,
Vol_{\eta} -\frac {6}{n}\int_M F\, \mathcal{C}(F)\vol,
\end{equation*}
where $u=F^2$, $(\nabla^2f)_{0}$ is the traceless part of horizontal Hessian of $f$ and $\mathcal{C}$ is the CR-Paneitz operator of $M$.
\end{prop}
We refer to Section \ref{s:CR} for the relevant notation and definitions. As a consequence of Proposition \ref{p:CR energy ineq} we recover the monotonicity of the entropy energy shown previously in \cite{CW10}. We note that one of my motivations to consider the problem was the application of the CR version of the monotonicity of the entropy like energy \cite[Lemma  3.3]{CW10} in obtaining (non-optimal) estimate on the bottom  of the $L^2$ spectrum of the CR sub-Laplacian. However, the proof of \cite[Corollary 1.9 and Section 6]{CW10} is not fully justified  since  \cite[Lemma3.3]{CW10} is proved for a compact manifold. It should be noted that a proof of S-Y Cheng's type (even non-optimal) estimate in a sub-Riemannian setting, such as CR or qc-manifold, is an interesting problem in particular because of  the lack of general comparison theorems.

 We conclude by mentioning another proof of the monotonicity of the energy in the recent preprint \cite{IP16}, which was the result of a past collaborative work with Ivanov  and Petkov. Remarkably, \cite{CW10} is also not acknowledged in \cite{IP16} despite the fact that the calculations in \cite{IP16} came after I introduced to Ivanov many of the interesting (sub-Riemannian) comparison problems and drew their attention to \cite{CW10}. While I can hardly wish to  be associated with \cite{IP16}, a quick look shows the line for line substantial overlap of  \cite[Section 3]{IP16} with Chang and Wu' proof \cite[Lemma 3.3]{CW10}, the publication of collaborative work without a discussion with all sides is notable.  Therefore, I decided to give my  independent approach to the problem.

{\bf Acknowledgments.}
Thanks are due to Qi Zhang for insightful comments on the S-Y Cheng eigenvalue estimate  during the	Beijing Workshop on Conformal Geometry and Geometric PDE in 2015, and to Jack Lee and Ben Chow for some useful discussions. The author also acknowledges  the  support of the Simons Foundation  grant \#279381.

\section{Proofs of the Propositions}
\subsection{Some preliminaries}\label{ss:prelim}
Throughout this section  $M$ will be a qc manifold  of dimension $4n+3$, \cite{Biq2}, with horizontal space $H$  locally given as the kernel of a 1-form $\eta=(\eta_1,\eta_2,\eta_3)$ with values in $\mathbb{R}^3$, and Biquard connection $\nabla$ with torsion $T$. Below we record some of the properties needed for this paper, see also  \cite{Biq1} and \cite{IV3} for a more expanded exposition.

The $Sp(n)Sp(1)$ structure on $H$ is
fixed by a positive definite symmetric tensor $g$ and a rank-three bundle $\mathbb{Q}$ of endomorphisms of $H$ locally generated by three almost complex structures $I_1,I_2,I_3$ on $H$ satisfying the identities of the imaginary unit quaternions and also the conditions $$g(I_s.,I_s.)=g(.,.)\qquad \text{and} \qquad 2g(I_sX,Y)\ =\
d\eta_s(X,Y).$$

Associated with the Biquard connection  is the vertical space $V$, which is complementary to $H$ in $TM$. In the case $n=1$ we shall make the usual assumption of existence of Reeb vector fields $\xi_1, \xi_2, \xi_3$, so that the connection is defined following D. Duchemin \cite{D}. The  fundamental 2-forms $\omega_s$ of the fixed qc
structure  will be denoted by $\omega_s$,
\begin{equation*}  
2\omega_{s|H}\ =\ \, d\eta_{s|H},\qquad
\xi\lrcorner\omega_s=0,\quad \xi\in V.
\end{equation*}

In order to give some idea of the involved quantities we list a few more essential for us details. Recall that
 $\nabla$ preserves the decomposition $H\oplus V$ and the $
Sp(n)Sp(1)$ structure on $H$, $$\nabla g=0, \quad \nabla\Gamma(\mathbb{Q})
\subset \Gamma( \mathbb{Q})$$ and its torsion on $H$ is given by
$T(X,Y)=-[X,Y]_{|V}$. Furthermore, for a vertical field $\xi\in V$, the endomorphism
$T_{\xi}\equiv T(\xi,.)_{|H}$ of $H$ belongs to the space $ (sp(n)\oplus sp(1))^{\bot}\subset
gl(4n)$ hence $ T({\xi}, X,Y)=g(T_{\xi}X,Y)$ is a well defined tensor field. The two $Sp(n)Sp(1)$-invariant trace-free symmetric
2-tensors $T^0(X,Y)= g((T_{\xi_1}^{0}I_1+T_{\xi_2}^{0}I_2+T_{
\xi_3}^{0}I_3)X,Y)$, $U(X,Y) =g(uX,Y)$ on $H$, introduced in
\cite{IMV},  satisfy
\begin{equation}  \label{propt}
\begin{aligned} T^0(X,Y)+T^0(I_1X,I_1Y)+T^0(I_2X,I_2Y)+T^0(I_3X,I_3Y)=0, \\
U(X,Y)=U(I_1X,I_1Y)=U(I_2X,I_2Y)=U(I_3X,I_3Y). \end{aligned}
\end{equation}
Note that when $n=1$, the tensor $U$ vanishes. The tensors $T^0$ and $U$ determine completely the torsion endomorphism  due to
the identity \cite[Proposition~2.3]{IV} $$%
4T^0(\xi_s,I_sX,Y)=T^0(X,Y)-T^0(I_sX,I_sY),$$ which in view of \eqref{propt} implies
\begin{equation}\label{need1}
\sum_{s=1}^3T(\xi_s,I_sX,Y)=
T^0(X,Y)-3U(X,Y).
\end{equation}
The curvature of the Biquard connection is $R=[\nabla,\nabla]-\nabla_{[\ ,\ ]}$ with  {qc-Ricci tensor} and \textit{normalized} qc-scalar
curvature, defined by respectively by
\begin{equation*}  \label{e:horizontal ricci}
Ric(X,Y)=\sum_{a=1}^{4n}{g(R(e_a,X)Y,e_a)},  \qquad
8n(n+2)S=\sum_{a=1}^{4n}Ric(e_a,e_a).
\end{equation*}
According to \cite{Biq1} the Ricci tensor restricted to $H$ is a
symmetric tensor. Remarkably, the torsion tensor determines  the qc-Ricci tensor of the Biquard connection on $M$ in view of the formula, \cite{IMV},
\begin{equation}\label{e:Ric by torsion}
Ric(X,Y) \ =\ (2n+2)T^0(X,Y) +(4n+10)U(X,Y)+\frac{S}{4n}g(X,Y).
\end{equation}
Finally, $\vol$ will denote the volume form
\begin{equation}\label{e:volume form}
\vol=\eta_1\wedge\eta_2\wedge\eta_3\wedge\Omega^n,
\end{equation}
where $\Omega=\omega_1\wedge\omega_1+\omega_2\wedge\omega_2+\omega_3\wedge\omega_3$ is the {fundamental 4-form}. We note the integration by parts formula
\begin{equation}  \label{div}
\int_M (\nabla^*\sigma)\,\, Vol_{\eta}\ =\ 0,
\end{equation}
where the (horizontal) divergence of a horizontal
vector field $\sigma\in\Lambda^1\, (H)$ is given by $\nabla^* \sigma =-tr|_{H}\nabla\sigma= -\nabla \sigma(e_a,e_a) $ for an orthonormal frame $\{e_a\}_{a=1}^{4n}$ of the horizontal space.

\subsection{Proof of Proposition \protect\ref{p:Paneitz estimate}}

We start with a formula for the change of the dependent
function in the $P$-function of $f$. To this effect, with $f=f(F)$, a short
calculation shows the next identity
\begin{multline*}
\nabla^3 f(Z,X,Y)=f^{\prime }\nabla^3F(Z,X,Y)+f^{\prime \prime \prime
}dF(Z)dF(X)dF(Y) \\
+f^{\prime \prime }\nabla^2F(Z,X)dF(Y) +f^{\prime \prime }\nabla^2F(Z,Y)dF(X)+f^{\prime
\prime }\nabla^2F(X,Y)dF(Z).
\end{multline*}
Recalling definition \eqref{e:P form} we obtain
\begin{multline}  \label{e:P from change variable}
P_f(Z)=f^{\prime }P_F(Z) + f^{\prime \prime \prime }|\Cr F|^2dF(Z) + 2 f^{\prime
\prime 2}F (Z, \nabla F)+ f^{\prime \prime }(\Delta F)dF(Z) \\
+f^{\prime \prime }\sum_{s=1}^3 g(\nabla^2 F, \omega_s)dF(I_s Z),
\end{multline}
which implies the  identity
\begin{equation}  \label{e:Paneitz change variable}
P_f(\nabla f)=f^{\prime 2}P_F(\nabla F) + f^{\prime }f^{\prime \prime \prime
 }|\Cr F|^4 + 2f^{\prime }f^{\prime \prime }\nabla ^2F (\nabla F, \nabla F)+f^{\prime
}f^{\prime \prime }|\Cr F|^2\Delta F.
\end{equation}

In our case, since we are interested in expressing the integral of $%
uP_f(\nabla f)=e^{-f}P_f(\nabla f)$ in terms of the integral of a $P$%
-function of some function, equation \eqref{e:Paneitz change variable} leads
to the ordinary differential equation $u\left (-\frac {u^{\prime }}{u}%
\right)^2=const$. Therefore, we let $u=F^2$ and find
\begin{equation}  \label{e:Paneitz change variable F}
uP_f(\nabla f)=4P_F(\nabla F) + 8F^{-2}|\nabla F|^4 -8F^{-1} \nabla^2F
(\nabla F, \nabla F)-4F^{-1}|\nabla F|^2\Delta F.
\end{equation}
Now, the last three terms will be expressed back in the variable $f$ which gives
\begin{equation}  \label{e:Paneitz change variable f}
uP_f(\nabla f)=4P_F(\nabla F) + \left [-\frac 14|\nabla f|^4+\frac 12
|\nabla f|^2\Delta f+ \nabla^2f (\nabla f, \nabla f)\right]u
\end{equation}
At this point, we integrate the above identity and then apply the (integration
by parts) divergence formula \eqref{div} in order to show
\begin{equation*}
\int_M \nabla^2f (\nabla f, \nabla f)u \, Vol_{\eta}=\frac 12\int_M \left[
|\nabla f|^4-|\nabla f|^2\Delta f\right]u \, Vol_{\eta},
\end{equation*}
which leads to \eqref{e:Paneitz estimate}. The proof of Proposition \ref{p:Paneitz estimate} is complete.

\subsection{Proof of Proposition \protect\ref{p:energy ineq}}

The initial step is identical to the Riemannian case \cite{Ni04}, so we skip the computations. Let $w=2\Delta f - |\nabla
f|^2$. Using the heat equation and integration by parts, exactly as in the Riemannian case, we have
\begin{equation}\label{e:E'}
\ddt \mathcal{E}(t)=\int_M(\partial_t-\Delta)(uw)\, Vol_{\eta}
\end{equation}
and also
\begin{equation}\label{e:key 1}
(\partial_t-\Delta)(uw)=\left [2g\left ( \nabla\left (\Delta f \right),
\nabla f\right)-\Delta|\nabla f|^2 \right]u.
\end{equation}
Next, we apply the qc Bochner formula \cite{IPV1, IPV2}
\begin{multline*}
\frac12\triangle |\nabla f|^2=|\nabla^2f|^2+g\left (\nabla (\triangle f),
\nabla f \right )+2(n+2)S|\nabla f|^2 \\
+2(n+2)T^0(\nabla f,\nabla f) +4(n+1)U(\nabla f,\nabla f) + 4R_f(\nabla f),
\end{multline*}
where
\begin{equation*}
R_f(Z)=\sum_{s=1}^3\nabla^2f(\xi_s,I_sZ).
\end{equation*}
Therefore,
\begin{multline}  \label{e:dt-lap}
\frac 12(\partial_t-\Delta)(uw) =\big[-|\nabla^2f|^2-2(n+2)S|\nabla
f|^2-2(n+2)T^0(\nabla f,\nabla f) \\
-4(n+1)U(\nabla f,\nabla f) - 4R_f (\nabla f)\big]u
\end{multline}
The next step is the computation of $\int_M R_f (\nabla f)u\, Vol_{\eta}$ in two ways as was done in \cite{IPV1,IPV2} for the Lichnerowicz type first eigenvalue lower bound but integrating with respect to $\vol$ rather than $u\vol$ as we need to do here. For ease of reading we will follow closely \cite[Section
8.1.1]{IV14} but notice the opposite convention of the sub-Laplacian in \cite[Section 8.1.1]{IV14}.
First with the help of the $P$-function, working similarly to \cite%
[Lemma 3.2]{IPV2} where the integration was with respect to $\vol$, we have
\begin{multline}  \label{e:Rf u1}
\int_M R_f (\nabla f)u\, Vol_{\eta}=\int_M [- \frac{1}{4n}P_n(\nabla f)-%
\frac{1}{4n}(\triangle f)^2-S|\nabla f|^2 \\
+ \frac{n+1}{n-1}U(\nabla f,\nabla f)]u\, Vol_{\eta} +\frac {1}{4n}\int_M
|\nabla f|^2(\Delta f)u\, Vol_{\eta},
\end{multline}
with the convention that in the case $n=1$ the formula is understood by formally dropping the term
involving (the vanishing) tensor $U$. Notice the appearance of a "new" term
in the last integral in comparison to the analogous formula in \cite[Section
8.1.1, p. 310]{IV14}. Indeed, taking into account the $Sp(n)Sp(1)$ invariance of $R_f(\nabla f)$
and Ricci's identities we have, cf. \cite[Lemma 3.2]{IPV2},
\begin{equation*}
R_f(X)=-\frac {1}{4n}\sum_{s=1}^3\sum_{a=1}^{4n} \nabla^3 f(I_sX, e_a, I_s e_a)+\left [ T^0(X,\nabla f)-3U(X, \nabla f)\right]
\end{equation*}
hence \eqref{e:P form} gives
\begin{equation*}
uR_f (\nabla f)= \big[- \frac{1}{4n}P_n(\nabla f)-S|\nabla f|^2 + \frac{n+1}{n-1}%
U(\nabla f,\nabla f)\big]u\newline
+\frac {1}{4n}\sum_{a=1}^{4n}\nabla^3f(\nabla f, e_a, e_a)u.
\end{equation*}
An integration by parts shows the validity of \eqref{e:Rf u1}.

On the other hand, we have
\begin{multline}  \label{e:Rf u2}
\int_MR_f(\nabla f) u\, Vol_{\eta} =-\int_M \left[\frac {1}{4n}\sum_{s=1}^3
g(\nabla^2 f, \omega_s)^2 +T^0(\nabla f,\nabla f)-3U(\nabla f,\nabla f)%
\right]u\vol,
\end{multline}
which other than using different volume forms is identical to the second formula in \cite[Section 8.1.1, p. 310]%
{IV14}.  Indeed, following \cite[Lemma 3.4]{IPV1}, using Ricci's identity \[\nabla^2f
(X,\xi_s)-\nabla^2f(\xi_s,X)=T(\xi_s,X,\nabla f)\] and \eqref{need1}, we have
\[
R_f(\Cr f)=\left(\sum_{s=1}^3\nabla^2f(I_s\Cr f, \xi_s) \right)-\left[ T^0(\Cr f, \Cr f)-3U(\Cr f, \Cr f)\right]
\]
An integration by parts gives \eqref{e:Rf u2}, noting
the term  $\sum_{s=1}^3
df(\xi_s)df(I_s\nabla f)=0$ and taking into account
that by Ricci's identity
$$\nabla^2f
(X,Y)-\nabla^2f(Y,X)=-2\sum_{s=1}^3\omega_s(X,Y)df(\xi_s)$$
we have
$
g(\nabla^2f , \omega_s) {=}\sum_{a=1}^{4n}\nabla^2f(e_a,I_se_a)=-4ndf(\xi_s).
$

Now, working as in \cite[Section 8.1.1, p. 310]{IV14}, we subtract \eqref{e:Rf u2}
and three times formula \eqref{e:Rf u1} from \eqref{e:dt-lap} which brings
us to the following identity
\begin{multline}  \label{e:E formula 1}
\frac 12 \frac{d}{dt}\mathcal{E}(t)=\int_M\big[ -|(\nabla^2f)_{0}|^2-\frac {2n+1%
}{2}\mathcal{L}(\nabla f,\nabla f)\big]u\, Vol_{\eta} \\
+\frac {1}{4n}\int_M \big[3P_f(\nabla f)+2(\Delta f)^2-3|\nabla f |^2\Delta
f \big]u \vol,
\end{multline}
where $|(\nabla^2f)_{0}|^2$ is the square of the norm of the traceless part of the horizontal Hessian
\begin{equation*}
|(\nabla^2 f)_0|^2=|\nabla^2f|^2-\frac{1}{4n}\Big[(\triangle
f)^2+\sum_{s=1}^{3}[g(\nabla^2f,\omega_s)]^2\Big].
\end{equation*}
Next, we consider $\int_M \big[2(\Delta f)^2-3|\nabla f |^2\Delta f \big]u\,
Vol_{\eta}$. Using the heat equation we have the identical to the Riemannian case relation
\begin{equation*}
\frac{d}{dt} \mathcal{E}(t)=\frac{d}{dt} \int_M w \Delta u\,
Vol_{\eta}=\int_M \big(-2(\Delta f)^2 + 3|\nabla f|^2\Delta f - |\nabla f|^4%
\big )u\, Vol_{\eta},
\end{equation*}
hence
\begin{equation}\label{e:E formula 2}
\int_M \big( 2(\Delta f)^2-3|\nabla f|^2\Delta f \big) u\, Vol_{\eta}= -%
\frac{d}{dt} \mathcal{E}(t)-\int_M |\nabla f|^4u\, Vol_{\eta}.
\end{equation}
A substitution of the above formula in \eqref{e:E formula 1} gives
\begin{multline*}
\frac {2n+1}{4n}\frac{d}{dt} \mathcal{E}(t)=\int_M\big[ -|(\nabla^2f)_0 |^2-%
\frac {2n+1}{2}\mathcal{L}(\nabla f,\nabla f)\big]u\, Vol_{\eta}
+\frac {1}{4n}\int_M \big[3P_f(\nabla f)- |\nabla f|^4\big]u\, Vol_{\eta}.
\end{multline*}
Finally, we invoke Proposition \ref{p:Paneitz estimate} in order to complete the
proof.

\section{The CR case}\label{s:CR}
In this section we prove the monotonicity formula in the CR case stated in Proposition \ref{p:CR energy ineq} following the method we employed in the qc case. This implies the monotonicity of the entropy like energy which was proved earlier in \cite{CW10}.

 Throughout the section $M$ will be a $(2n+1)$-dimensional strictly pseudoconvex (integrable) CR manifold with a fixed pseudohermitian structure defined by a contact form $\eta$ and complex structure $J$ on the horizontal space $H=Ker \, \eta$. The fundamental 2-form is defined by  $\omega=\frac 12 \eta$ and the  Webster metric is $g(X,Y)=-\omega (JX,Y)$ which is extended to a Riemannian metric on $M$ by declaring that the Reeb vector field associated to $\eta$ is of length one and orthonormal to the horizontal space. We shall denote by $\Cr$ the associated Tanaka-Webster connection \cite{T} and \cite{W, W1}, while $\triangle u=tr^g(\nabla^2 u)$ will be the negative sub-Laplacian with the trace taken with respect to an orthonormal basis of the horizontal $2n$-dimensional space. Finally, we define the Ricci type tensor
 \begin{equation}\label{e:Lichnerowicz tensor CR}
 \mathcal{L}(X,Y)=\rho(JX,Y)+2nA(JX,Y)
 \end{equation}
 recalling that on a CR manifold we have
 \begin{equation}  \label{rid}
Ric(X,Y)=\rho(JX,Y)+2(n-1)A(JX,Y),
\end{equation}
where $\rho$ is the $(1,1)$-part of the pseudohermitian Ricci
tensor (the Webster Ricci tensor) while the
$(2,0)+(0,2)$-part is the Webster torsion $A$, see \cite[Chapter 7]{IV3} for the  expressions in real coordinates of these known formulas \cite{W,W1}, see also \cite{DT}.
 
 With the above convention in place, as in \cite{CW10}, for a positive solution of \eqref{e:heat}  we consider the entropy \eqref{e:N definition} and energy \eqref{enf}, where $\vol=\eta\wedge(d\eta)^{2n}$. 
 
We turn to the proof of Proposition \ref{p:CR energy ineq}. For a function $f$ we define the one form,
\begin{equation}  \label{e:Pdef}
 P_{f}(X)=\nabla ^{3}f(X,e_{b},e_{b})+\nabla
^{3}f(JX,e_{b},Je_{b})+4nA(X,J\nabla f)
\end{equation}%
so that the  fourth order  CR-Paneitz operator  is given by
\begin{multline}  \label{e:Cdef}
C(f)=-\nabla ^{\ast }P=(\nabla_{e_a} P)({e_a})=\nabla ^{4}f(e_a,e_a,e_{b},e_{b})+\nabla
^{4}f(e_a,Je_a,e_{b},Je_{b})\\
-4n\nabla^* A(J\nabla f)-4n\,g(\nabla^2 f,JA).
\end{multline}
By \cite{GL88}, when $n>1$ a  function $f\in \mathcal{C}^3(M)$ satisfies the equation $Cf=0$ iff $f$ is CR-pluriharmonic. Furthermore, the CR-Paneitz operator is non-negative,
\[
\int_M f\cdot Cf \vol=-\int_MP_f(\gr) \vol\geq 0.
\]
 On the other hand, in the three dimensional case the positivity condition is a CR invariant since it is independent of the choice of the contact form by the conformal invariance of $C$ proven in \cite{Hi93}.
 
We turn to the proof of Proposition \ref{p:CR energy ineq}. Taking into account \eqref{e:key 1} and the CR Bochner formula \cite{Gr}, 
\begin{multline}\label{e:bohh}
\frac12\triangle |\nabla f|^2=|\nabla^2 f|^2+g(\nabla(\triangle f),\nabla f)+Ric(\nabla
f,\nabla f)+2A(J\nabla f,\nabla f) 
+ 4R_f(\Cr f),
\end{multline}
where $R_f(Z)=\nabla df(\xi,JZ)$,  see \cite[Section 7.1]{IV14} and references therein but note the opposite sign of the sub-Laplacian, we obtain the next identity
\begin{equation}  \label{e:dt-CR lap}
\frac 12(\partial_t-\Delta)(uw) =\big[-|\nabla^2f|^2-Ric(\nabla f,\nabla f)
-2A(\nabla f,\Cr \nabla f) - 4R_f (\nabla f)\big]u.
\end{equation}
Since \eqref{e:E'} still holds, working as in the qc case we compute $\int_M R_F(\Cr f)u\vol$ in two ways \cite[Lemma 4]{Gr} and \cite[Lemma 8.7]{IVO}  following the exposition \cite{IV14}.

From Ricci's identity  $$\nabla^2f (X,Y)-\nabla^2f(Y,X)=-2\omega(X,Y)df(\xi)$$ it follows $df(\xi)=-\frac {1}{2n}g(\nabla^2 f, \omega)$. Hence 
\[
\nabla^2f(JZ, \xi)=-\frac {1}{2n}\sum_{b=1}^{2n}\nabla^3 f(JZ, e_b,Je_b),
\]
where $\{e_b\}_{b=1}^{2n}$ is an orthonormal basis of the horizontal space. Applying Ricci's identity
\[
\nabla^2f (X,\xi)-\nabla^2f(\xi,X)=A(X,\nabla f)
\]
it follow
\begin{equation}\label{e:CR rem}
R_f(Z)=\nabla ^{2}f(\xi ,JZ)=-\frac{1}{2n}\sum_{b=1}^{2n}\nabla^3 f(JZ, e_b,Je_b) -A(JZ,\nabla f).
\end{equation}
Taking into account \eqref{e:Pdef} the last formula gives 
\[
R_f(Z)=-\frac{1}{2n}P_f(Z) +A(JZ,\nabla f)+\frac{1}{2n}\sum_{b=1}^{2n}\nabla^3 f(Z, e_b,e_b).
\]
Now, an integration by parts shows the next identity
\begin{multline}\label{e:Rf u1 CR}
\int_M R_f(\Cr f) u\vol=\int_M\big[-\frac{1}{2n}P_f(\Cr f) +A(J\nabla f,\nabla f)-\frac {1}{2n}(\Delta f)^2  +\frac {1}{2n}|\Cr f|^2 (\Delta f) \big]u\vol.
\end{multline}
On the other hand, using again \eqref{e:CR rem} but now we integrate and then use integration by parts we have
\begin{equation}\label{e:Rf u2 CR}
\int_M R_f(\Cr f) u\vol=\int_M\big[ -\frac {1}{2n}g(\nabla^2 f,\omega)^2-A(J\nabla f, \nabla f)\big ] u\vol.
\end{equation}

At this point, exactly as in the qc case, we subtract \eqref{e:Rf u2 CR}
and three times formula \eqref{e:Rf u1 CR} from \eqref{e:dt-CR lap}, which gives
\begin{equation*}
\mathcal{E}'(t)=-\int_M\big[ |(\nabla^2f)_{0}|^2 +\mathcal{L}(\nabla f,\nabla f)\big]u\vol \\
+\frac {1}{2n}\int_M \big[3P_f(\nabla f)+2(\Delta f)^2-3|\nabla f |^2\Delta
f \big]u\vol,
\end{equation*}
where  $|(\nabla^2f)_{0}|^2$ is the square of the norm of the traceless part of the horizontal Hessian
\begin{equation*}
|(\nabla^2 f)_0|^2=|\nabla^2f|^2-\frac{1}{2n}\Big[(\triangle
f)^2+g(\nabla^2f,\omega)^2\Big].
\end{equation*}
Taking into account that the formulas in Proposition \ref{p:Paneitz estimate} and \eqref{e:E formula 2} hold unchanged we complete the proof.

\end{document}